\documentclass[10pt,a4paper]{amsart}

\input xy
\xyoption{all}

\newtheorem{thm}{Theorem}[section]
\newtheorem{lemma}[thm]{Lemma}
\newtheorem{cor}[thm]{Corollary}
\newtheorem{pro}[thm]{Proposition}
\newtheorem{example}[thm]{Example}
\newtheorem{remark}[thm]{Remark}
\newtheorem{ddef}[thm]{Definition}

\def\ext{\mathop{\rm Ext}\nolimits}

\def\tor{\mathop{\rm Tor}\nolimits}
\def\hom{\mathop{\rm Hom}\nolimits}

\def\a{{\alpha}}
\def\b{{\beta}}

\def\homgr{\mathop{\rm Homgr}\nolimits}

\def\C{{\mathcal C}}

\def\om{{\omega}}

\def\s{{\sigma}}
\def\a{{\alpha}}

\def\b{{\beta}}

\def\d{{\delta}}

\def\iq{{\mathfrak q}}

\def\ip{{\mathfrak p}}

\def\ra{{\rightarrow}}

\def\fini{{$\quad\quad\Box$}}

\newcommand{\bd}{\begin{ddef}}
\newcommand{\ed}{\end{ddef}}
\newcommand{\bt}{\begin{thm}}
\newcommand{\et}{\end{thm}}
\newcommand{\bl}{\begin{lemma}}
\newcommand{\el}{\end{lemma}}
\newcommand{\bco}{\begin{cor}}
\newcommand{\eco}{\end{cor}}
\newcommand{\bp}{\begin{pro}}
\newcommand{\ep}{\end{pro}}
\newcommand{\bex}{\begin{example}}
\newcommand{\eex}{\end{example}}
\newcommand{\brm}{\begin{remark}}
\newcommand{\erm}{\end{remark}}
\newcommand{\bconj}{\begin{conj}}
\newcommand{\econj}{\end{conj}}

\newcommand{\beqn}{\begin{eqnarray*}}
\newcommand{\eeqn}{\end{eqnarray*}}
\newcommand{\beq}{\begin{eqnarray}}
\newcommand{\eeq}{\end{eqnarray}}
\newcommand{\been}{\begin{enumerate}}
\newcommand{\eeen}{\end{enumerate}}

\begin{document}

\title{A duality theorem for generalized local cohomology}

\author[Chardin, Divaani-Aazar]{Marc Chardin and Kamran Divaani-Aazar}

\address{M. Chardin, Institut Math\'ematiques de Jussieu,
175 rue du Chevaleret,
F-75013, Paris, France.}
\email{chardin@math.jussieu.fr}

\address{K. Divaani-Aazar, Department of Mathematics, Az-Zahra
University, Vanak, Post Code 19834, Tehran, Iran.} \email{kdivaani@ipm.ir}



\begin{abstract} We prove a duality theorem for graded algebras over a field that implies several 
known duality results : graded local duality, versions of Serre duality for local cohomology and of Suzuki duality for generalized local cohomology, and Herzog-Rahimi bigraded duality. 
\end{abstract}

\maketitle

\section{Introduction}

Several duality results for local cohomology (notably Grothendieck and Serre dualities) and 
for generalized local cohomology (Suzuki duality) are classical results commonly used in
commutative algebra. A more recent duality theorem was established by Herzog and Rahimi in the
context of bigraded algebras. We here present a duality result that encapsulates all these 
duality statements. Before stating this theorem, we need  to fix some notation.

Let $\Gamma$ be an abelian group and $R$ be a graded polynomial ring over a commutative Noetherian ring $R_0$. 

For a graded $R$-ideal $I$, we recall from \cite{CCHS} that the grading of $R$ is $I$-sharp if
$H^i_I(R)_\gamma$ is a finitely generated $R_0$-module for any $\gamma \in \Gamma$ and any $i$. If this is the case,
then $H^i_I(M)_\gamma$ is also a finitely generated $R_0$-module for any finitely generated graded 
$R$-module $M$, any $\gamma \in \Gamma$ and any $i$.

In the case $S=S_0 [x_1,\ldots ,x_m,y_1,\ldots y_n]$ is a graded quotient of  $R=R_0 [X_1,\ldots ,X_m,Y_1,\ldots Y_n]$,  $\ip := (x_1,\ldots ,x_m)$ and $\iq := (y_1,\ldots y_n)$, we define the  grading of $S$ to be $\ip$-sharp (equivalently $\iq$-sharp) if the following condition is satisfied :
for all $\gamma \in \Gamma$, 
$$
| \{ (\alpha , \beta )\in {\bf Z}_{\geq 0}^{m}\times {\bf Z}_{\geq 0}^{n}\ :\quad \sum_i \a_i\deg (x_i )=\gamma +\sum_j \b_j \deg (y_j )\}| <\infty .
$$
This condition on the degrees of the variables is equivalent to the fact that the grading of
$R$ is $(X_1,\ldots ,X_m)$-sharp (equivalently $(Y_1,\ldots ,Y_n)$-sharp) by \cite[1.2]{CCHS}.

We set ${^{*}\hom}_{R}(M,R_0)$ for the graded $R_0$-dual of a  graded $R$-module $M$. This module is graded and ${^{*}\hom}_{R}(M,R_0)_\gamma =\hom_{R_0}(M_{-\gamma},R_0)$. 

The notion of generalized local cohomology was given by Herzog in his
Habilitationsschrift \cite{H}. For an ideal $I$ in a commutative Noetherian ring
$R$ and $M$ and $N$ two $R$-modules, the $i$-th generalized local
cohomology module of $M$ and $N$ with respect to $I$ is denoted by $H^i_{I}(M,N)$. 
One has a natural isomorphism $H^i_{I}(R,N)\simeq H^i_{I}(N)$ and a spectral sequence
$H^i_I (\ext^j_R(M,N))\Rightarrow H^{i+j}_I(M,N)$.

Our main result is the following duality result,
\bt
Let $S=k[x_1,\ldots ,x_m,y_1,\ldots y_n]$ be a Cohen-Macaulay $\Gamma$-graded algebra over a field $k$. Let
$\ip := (x_1,\ldots ,x_m)$ and $\iq := (y_1,\ldots y_n)$. Assume that the grading of $S$ is $\ip$-sharp (equivalently $\iq$-sharp),
and $M$ and $N$ are finitely generated graded $S$-modules such that either $M$ has finite projective dimension  or $N$ 
has finite projective dimension  and $\tor^S_i (M,\om_S)=0$ for $i>0$, then
$$
 H^{i}_{\ip}(M,N)\simeq {^{*}\hom}_{S}(H^{\dim S-i}_{\iq}(N,M\otimes_S \om_S),k).
 $$
 \et

Several previously known duality results are particular cases of this theorem. Namely,

(i) graded local duality in the case that the base ring is a field corresponds to the case where $\iq =0$ and $M=S$,

(ii) Serre duality follows from the case where $\ip =0$ and $M=S$, {\it via} the spectral sequence
$H^i_\iq (\ext^j_S(N,\om_S))\Rightarrow H^{i+j}_\iq (N,\om_S)$,

(iii) Suzuki duality in the context of graded algebras over a field (see \cite[3.1]{CD}) corresponds to the case where $\iq =0$,

(iv) Herzog-Rahimi spectral sequence corresponds to the case where $M=S$ (see \cite{HR} for the
standrard bigraded case, and \cite{CCHS} for the general case) using the spectral sequence
$H^i_\iq (\ext^j_S(N,\om_S))\Rightarrow H^{i+j}_\iq (N,\om_S)$.

\bigskip

{\it This research took place when the second named author was on sabbatical leave. 
He would like to express his deep thanks to Institut Math\'ematiques de Jussieu for its kind hospitality.}

\section{The duality theorem}

If $F_\bullet$ and
$G_\bullet$ are two complexes of $R$-modules, $\homgr_{R}(F_{\bullet}, G_{\bullet})$ is the cohomological complex
with modules $C^i=\prod_{p-q=i}\hom_R(F_p,G_q)$. If either $F_{\bullet}$ or $G_{\bullet}$ is finite, then
$\homgr_{R}(F_{\bullet}, G_{\bullet})$ is the totalisation of the double complex with $C^{p,-q}=\hom_R(F_p,G_q)$.

We first consider the case where $R:=k[X_1,\ldots ,X_m,Y_1,\ldots ,Y_n]$ is a polynomial ring over a field $k$, 
$\ip :=(X_1,\ldots ,X_m)$ and $\iq :=(Y_1,\ldots ,Y_n)$. 

From this point on, we will assume that $R$ is $\Gamma$-graded, for an abelian group $\Gamma$, and that
the grading of $R$ is $\ip$-sharp (equivalently $\iq$-sharp).
One has $\om_R =R[-\s ]$, where $\s :=\sum_i \deg (X_i)+\sum_j \deg (Y_j)$.  

Let $\C^\bullet_\ip (M)$ be the \v Cech complex
$0\ra M\ra \oplus_i M_{X_i}\ra \cdots \ra M_{X_1\cdots X_m}\ra 0$ and  $\C^\bullet_\iq (M)$ be the \v Cech complex
$0\ra M\ra \oplus_i M_{Y_i}\ra \cdots \ra M_{Y_1\cdots Y_n}\ra 0$.  Denote by $C^{\bullet}\{ i\}$ the complex 
$C^{\bullet}$ shifted in cohomological degree by $i$ and set $\hbox{---}^*:={^{*}\hom}_{R}(\hbox{---},k)$. 
Recall that $\hbox{---}^*$ is exact on the category of graded $R$-modules. By  \cite[1.4 (b)]{CCHS} one has 

\bl\label{basicduality} If $F_{\bullet}$ is a graded complex of finite free $R$-modules, then 

(i) $H^{i}_{\ip}(F_{\bullet})=0$ for $i\not= m$,  $H^{j}_{\iq}(F_{\bullet})=0$ for $j\not= n$, 

(ii) there is a natural graded map of complexes $d: \C^\bullet_\ip (R)\{ m\}\ra \C^\bullet_\iq (\om_R)\{ n\}^*$ that induces a map of double complexes $\d : \C^\bullet_\ip (F_{\bullet})\{ m\}\ra  \C^\bullet_\iq (\hom_{R}(F_{\bullet},\om_R))\{ n\} ^*$ which gives rise to a functorial isomorphism
$$
H^{m}_{\ip}(F_{\bullet})\simeq H^{n}_{\iq}(\hom_{R}(F_{\bullet},\om_R))^*.
$$
\el

{\it Proof.} Remark that the map defined in \cite[1.4 (b)]{CCHS} is induced by the map from $\C^m_\ip (R)$ 
to $\C^n_\iq (\om_R)^*$ defined by $d(X^{-s-1}Y^{p})(X^{t}Y^{-q-1})=1$ if $s=t$ and $p=q$ and $0$ else (here $p$, $q$, $r$ and
$s$ are tuples of integers, with $p\geq 0$ and $t\geq 0$). Notice that $\C^\bullet_\ip (R)\{ m\}$ is zero in positive cohomological  degrees and $\C^\bullet_\iq (\om_R )\{ n\}^*$  is zero in negative cohomological  degrees. Hence this map extends uniquely to a map of complexes. \fini

\bl\label{functduality} Let $M$, $N$ and $N'$ be finitely generated graded $R$-modules and 
$f :N\ra N'$ be a graded homomorphism. For any integer $i$, there exist functorial isomorphisms $\psi_{N}^{i}$ and $\psi_{N'}^{i}$ 
that give rise to  a commutative diagram
$$
\small{
\xymatrix@!C@C=-58pt@H=20pt{
&H^{i-m}(H^{m}_{\ip}(\homgr_{R}(F_{\bullet}^{M}, F_{\bullet}^{N})))\ar'[d][dd]\ar_(.45){\simeq}^(.45){H^{i-m}(\d_{N})}[rr]&&\!H^{i-m}(H^{n}_{\iq}(\homgr_{R}(F_{\bullet}^{N}, F_{\bullet}^{M}\otimes_R \om_R))^{*})\ar[dd]\\
H^{i}_{\ip}(M,N)\ar^(.4){\psi_{N}^{i}}[rr]\ar[dd]\ar^{\simeq}[ur]&&H^{m+n-i}_{\iq}(N,M\otimes_R \om_R)^{*}\ar[dd]\ar^{\simeq}[ur]&\\
&H^{i-m}(H^{m}_{\ip}(\homgr_{R}(F_{\bullet}^{M}, F_{\bullet}^{N'})))\ar'[r]_(.8){\simeq}^(.8){H^{i-m}(\d_{N'})}[rr]&&\! H^{i-m}(H^{n}_{\iq}(\homgr_{R}(F_{\bullet}^{N'}, F_{\bullet}^{M}\otimes_R \om_R))^{*})\\
H^{i}_{\ip}(M,N')\ar^(.4){\psi_{N'}^{i}}[rr]\ar^{\simeq}[ur]&&H^{m+n-i}_{\iq}(N',M\otimes_R \om_R)^{*}\ar^{\simeq}[ur]&\\
}
}
$$
where the vertical maps are the natural ones.
\el

{\it Proof.} Let $F_{\bullet}^{M}$ be a minimal graded free $R$-resolution of $M$ and $f_{\bullet}: F_{\bullet}^{N}
\ra  F_{\bullet}^{N'}$ be a lifting of $f$ to minimal graded free $R$-resolutions of $N$ and $N'$. 

Recall that there is a spectral sequence $H^i(H^j_\ip (\homgr_R(F_{\bullet}^{M},F_{\bullet}^{N}))\Rightarrow H^{i+j}_\ip (M,N)$, by \cite[2.12]{CD} and \cite[\S 6, Th\'eor\`eme 1, b]{B}. 
As $H^{i}_{\ip}(F_{\bullet})=0$ for $i\not= m$ if $F_\bullet$ is a complex of free $R$-modules,  there is a
natural commutative diagram
$$
\xymatrix{
H^{i}_{\ip}(M,N)\ar^(.3){\simeq}[r]\ar^{can}[d]&H^{i-m}(H^{m}_{\ip}(\homgr_{R}(F_{\bullet}^{M}, F_{\bullet}^{N})))\ar^{can}[d]\\
H^{i}_{\ip}(M,N')\ar^(.3){\simeq}[r]&H^{i-m}(H^{m}_{\ip}(\homgr_{R}(F_{\bullet}^{M}, F_{\bullet}^{N'}))). \\
}
$$
As $F_{\bullet}^{M}$ is bounded, Lemma \ref{basicduality} provides the functorial isomorphisms
$$
\begin{array}{rcl}
H^{m}_{\ip}(\homgr_{R}(F_{\bullet}^{M}, F_{\bullet}^{N}))&\simeq &H^{n}_{\iq}(\hom_R(\homgr_{R}(F_{\bullet}^{M}, F_{\bullet}^{N}),\om_R))^*\\
&\simeq &H^{n}_{\iq}(\homgr_{R}(F_{\bullet}^{N}, F_{\bullet}^{M}\otimes_R \om_R))^*\\
\end{array}
$$
and equivalent ones for $N'$. This provides a commutative diagram
$$
\xymatrix{
H^{i}_{\ip}(M,N)\ar^(.25){\simeq}[r]\ar^{can}[d]&H^{i-m}(H^{n}_{\iq}(\homgr_{R}(F_{\bullet}^{N}, F_{\bullet}^{M}\otimes_R \om_R))^*)\ar^{can}[d]\\
H^{i}_{\ip}(M,N')\ar^(.25){\simeq}[r]&H^{i-m}(H^{n}_{\iq}(\homgr_{R}(F_{\bullet}^{N'}, F_{\bullet}^{M}\otimes_R \om_R))^*). \\
}
$$
As $k$ is a field, we have a natural isomorphism
$$
H^{i-m}(H^{n}_{\iq}(\homgr_{R}(F_{\bullet}^{N}, F_{\bullet}^{M}\otimes_R \om_R))^*)
\simeq H^{m-i}(H^{n}_{\iq}(\homgr_{R}(F_{\bullet}^{N}, F_{\bullet}^{M}\otimes_R \om_R))^*
$$
for $N$ and a similar one for $N'$. As $F_{\bullet}^{M}\otimes_R \om_R$ is acyclic, it provides a free $R$-resolution of $M\otimes_R \om_R$ and  we get a commutative diagram
$$
\small{
\xymatrix{
H^{m-i}(H^{n}_{\iq}(\homgr_{R}(F_{\bullet}^{N}, F_{\bullet}^{M}\otimes_R \om_R)))\ar^(.6){\simeq}[r]&H^{m+n-i}_{\iq}(N,M\otimes_R \om_R)\\
H^{m-i}(H^{n}_{\iq}(\homgr_{R}(F_{\bullet}^{N'}, F_{\bullet}^{M}\otimes_R \om_R)))\ar^{can}[u]\ar^(.6){\simeq}[r]&H^{m+n-i}_{\iq}(N',M\otimes_R \om_R)\ar^{can}[u],\\}
}
$$
that follows from same argument as in the first part of this proof, applied to $\iq ,N,M\otimes_R \om_R$ and $\iq ,N',M\otimes_R \om_R$ in place of $\ip ,M,N$ and $\ip ,M,N'$ .
This finally gives the claimed diagram of natural maps. \fini

\brm
Lemma \ref{functduality} holds more generally if we replace $k$ by a Gorenstein ring of dimension 0, and
assume that either $M$ has finite projective dimension or $N$ has finite projective dimension and $\tor_i^S(M,\om_S)=0$ for $i>0$. 
\erm

\bco\label{DFFCM}
Assume $S$ is a Cohen-Macaulay ring and a finitely generated graded $R$-module. Let  $F_\bullet$ be a graded 
complex  of finite free $S$-modules, and consider a graded  free $R$-resolution $F_{\bullet\bullet}$ of 
 $F_\bullet$. Then the map given by Lemma \ref{basicduality}
$$
\d :  \C^\bullet_\ip (F_{\bullet\bullet})\{ m\}\ra \C^\bullet_\iq (\hom_R(F_{\bullet\bullet},\om_R))\{ n\} ^*
$$ 
induces  a graded functorial map of double complexes 
$$
\C^{\bullet}_{\ip}(F_{\bullet})\{ \dim S\}\ra\C^{\bullet}_{\iq}(\hom_{S}(F_{\bullet},\om_S)) ^*
$$
which gives rise to isomorphisms
$$
H^{i}_{\ip}(F_{\bullet})\simeq {^{*}\hom}_{S}(H^{\dim S-i}_{\iq}(\hom_{S}(F_{\bullet},\om_S)),k).
$$
\eco

{\it Proof.} First notice that $S$ is a Cohen-Macaulay $R$-module and a graded $k$-algebra. By Lemma \ref{functduality}, the map $\d$ induces commutative diagrams
$$
\xymatrix{
H^{i}_{\ip}(R,F_{j-1})\ar^(.4){\simeq}[r]\ar[d]& H^{m+n-i}_{\iq}(F_{j-1},\om_R)^*\ar[d]\\
H^{i}_{\ip}(R,F_{j})\ar^(.4){\simeq}[r]\ar[d]&H^{m+n-i}_{\iq}(F_{j},\om_R)^*\ar[d]\\
H^{i}_{\ip}(R,F_{j+1})\ar^(.4){\simeq}[r]&H^{m+n-i}_{\iq}(F_{j+1},\om_R)^*\\
}
$$

This gives the claim together with the  isomorphisms $H^{i}_{\ip}(R,F_{\bullet})\simeq H^{i}_{\ip}(F_{\bullet})$ and $H^{m+n-i}_{\iq}(F_{\bullet},\om_R)\simeq H^{\dim S-i}_\iq (\hom_S (F_{\bullet},\om_S))$, where the second one follows from the Cohen-Macaulayness of $S$ {\it via} the collapsing
spectral sequence $H^i_\iq (\ext^j_R(F_{\bullet},\om_R ))\Rightarrow H^{i+j}_\iq (F_{\bullet},\om_R )$.\fini

\bl\label{tensext}
Assume $S$ is a Cohen-Macaulay ring and a finitely generated graded $R$-module. 
Let $M$ and $N$ be finitely generated 
graded $S$-modules  such that  $N$ has finite projective dimension and $\tor_i^S(M,\om_S)=0$ for $i>0$. Then the map $f\mapsto f\otimes 1_{\om_S}$ induces natural isomorphisms
$$
\ext^{i}_S(M,N)\simeq \ext^{i}_S(M\otimes_S \om_S,N\otimes_S \om_S).
$$
It follows that for any $S$-ideal $I$, $H^{i}_I (M,N)\simeq H^{i}_I(M\otimes_S \om_S,N\otimes_S \om_S)$ for all $i$.
\el

{\it Proof.}  The case where $M=N=S$ is a classical property of $\om_S$.  
If $M=S$ the result follows by induction on the projective dimension of
$N$, using the five Lemma and the exact sequences of Ext modules derived from the exact sequences $E : 0\ra \hbox{Syz}_1^S(N)\ra F\ra N\ra 0$ and $E\otimes_S \om_S$. Notice $E\otimes_S \om_S$ is exact since $\tor_i^S(P,\om_S )=0$ for $i>0$ if $P$ has finite projective dimension (see \cite[3.17, 4.5]{S}). Hence : $(*)\ \  \hom_S(\om_S ,N\otimes_S \om_S)\simeq N\ \hbox{and}\ \ext^i_S(\om_S ,N\otimes_S \om_S)=0,\ \forall i>0$.

For the general case, let $F_\bullet^M$ be a free $S$-resolution of $M$,  $F_\bullet^N$ be a finite free $S$-resolution of $N$ and $I^\bullet_{\om_S}$ be an injective $S$-resolution of
$\om_S$. The complex $F_\bullet^M\otimes_S \om_S$ is acyclic and, as $N$ has finite projective dimension, $F_\bullet^N\otimes_S I^{\bullet}_{\om_S}$ is 
a complex of injective $S$-modules whose only non zero homology is $N\otimes_S \om_S$ sitting in degree $0$, by \cite[\S 6, Th\'eor\`eme 1, a]{B}. It follows that 
$$
\begin{array}{rl}
\ext^{i}_S(M\otimes_S \om_S, N\otimes_S \om_S)
&\simeq H^i(\homgr_S(F_\bullet^M\otimes_S \om_S,F_\bullet^N\otimes_S I^{\bullet}_{\om_S}))\quad \hbox{by \cite[\S 6, Th\'eor\`eme 1, b]{B}},\\
&\simeq H^i(\homgr_S(F_\bullet^M,\hom_S (\om_S,F_\bullet^N\otimes_SI^\bullet_{\om_S}))\quad \hbox{by adjointness},\\
&\simeq H^i(\hom_S(F_\bullet^M,N)\quad \hbox{by {{\it ibid.}} and $(*)$},\\
&\simeq  \ext^{i}_S(M,N).\\
\end{array}
$$

It follows that the map from $\C^\bullet_I (\homgr_S(F_\bullet^M,F_\bullet^N))$ to 
$\C^\bullet_I (\homgr_S(F_\bullet^M\otimes_S \om_S,F_\bullet^N\otimes_S I^{\bullet}_{\om_S}
))$ induced by 
tensoring by $1_{\om_S}$, gives a morphism of spectral sequences which is an isomorphism
between  the first terms, which are respectively  $\C^p_I (\ext^q_S(M,N))$ and $\C^p_I (\ext^q_S(M\otimes_S \om_S,N\otimes_S \om_S))$.  Hence the total homology of these double complexes are isomorphic. This proves our claim by \cite[2.12]{CD}.
\fini

\bt
Let $S=k[x_1,\ldots ,x_m,y_1,\ldots y_n]$ be a Cohen-Macaulay $\Gamma$-graded algebra over a field $k$. Let
$\ip := (x_1,\ldots ,x_m)$ and $\iq := (y_1,\ldots y_n)$. Assume that the grading of $S$ is $\ip$-sharp (equivalently $\iq$-sharp),
and $M$ and $N$ are finitely generated graded $S$-modules such that either $M$ has finite projective dimension  or $N$ 
has finite projective dimension and $\tor^S_i (M,\om_S)=0$ for $i>0$, then
$$
 H^{i}_{\ip}(M,N)\simeq {^{*}\hom}_{S}(H^{\dim S-i}_{\iq}(N,M\otimes_S \om_S),k).
 $$
 \et

 {\it Proof.} Since the grading of $S$ is $\ip$-sharp, it follows that $H^{i}_{\ip}(M,N)\simeq H^{i}_{\ip}(M,N)^{**}$, for all $i$. Hence by Lemma \ref{tensext} it suffices to treat the case where $M$ has finite projective dimension. Let $F_{\bullet}^{M}$ and $F_{\bullet}^{N}$ be a minimal graded free $S$-resolutions of $M$ and $N$, respectively. Set $F_{\bullet}:=\homgr_{S}(F_{\bullet}^{M}, F_{\bullet}^{N})$. 
 
By Corollary \ref{DFFCM}, if $F_{\bullet\bullet}$ is a graded free $R$-resolution of $F_\bullet$, the map 
$$
\d :  \C^\bullet_\ip (F_{\bullet\bullet})\{ m\}\ra \C^\bullet_\iq (\hom_R(F_{\bullet\bullet},\om_R))\{ n\} ^*
$$
given by Lemma \ref{basicduality}
induces  a graded morphism of double complexes
$$
\C^{\bullet}_{\ip}(\homgr_{S}(F_{\bullet}^{M}, F_{\bullet}^{N}))\{ \dim S\}\ra\C^{\bullet}_{\iq}(\homgr_{S}(F_{\bullet}^{N}, F_{\bullet}^{M}\otimes_S \om_S)) ^*
$$
that is an isomorphism on the $E_1$ terms :
$$
H^{i}_{\ip}(\homgr_{S}(F_{\bullet}^{M}, F_{\bullet}^{N}))\simeq H^{\dim S-i}_{\iq}(\homgr_{S}(F_{\bullet}^{N}, F_{\bullet}^{M}\otimes_S \om_S))^*.
$$
It follows that these two spectral sequences have isomorphic terms at any level, in particular the same abutment. The first one abuts to $H^i_\ip (M,N)$ by \cite[2.12]{CD}.  The second one abuts to
$H^{\dim S-i}_{\iq}(N, M\otimes_S \om_S)^*$ by  \cite[\S 6, Th\'eor\`eme 1, b]{B} and \cite[2.12]{CD}, as $\tor^S_i (M,\om_S)=0$ for $i>0$ (because $M$ has finite projective dimension).\fini

\end{document}